\documentclass[12pt,a4paper]{amsart}

\usepackage{latexsym} 
 
\usepackage{graphicx}
\usepackage{epsfig}
\usepackage{mathrsfs}
\usepackage{amssymb}
\usepackage{amsthm}
\usepackage{amsfonts}
\usepackage{amsmath}
\usepackage{amstext}
\usepackage{amscd}
\usepackage{tikz}
\usepackage{epic}
\usepackage{eepic}
\usepackage{pstricks,pst-plot}

\numberwithin{equation}{section}

\theoremstyle{plain}
\newtheorem{thm}{Theorem}[section] 
\newtheorem{prop}[thm]{Proposition}
\newtheorem{cor}[thm]{Corollary}
\newtheorem{lem}[thm]{Lemma}

\newtheorem{theorem*}{Theorem}[]

\theoremstyle{definition}
\newtheorem{defn}[thm]{Definition}

\newtheorem{example}[thm]{Example}

\theoremstyle{remark}
\newtheorem{rem}[thm]{Remark}
\theoremstyle{property}

\newcommand{\N}{\mathbb{N}}
\newcommand{\R}{\mathbb{R}}

\DeclareMathOperator{\grad}{grad\,}

\def\accentclass@{7}
\def\makeacc@#1#2{\def#1{\mathaccent"\accentclass@#2 }}
\makeacc@\cir{017}

\title{The Lipschitz type of the Geometric Directional Bundle}

\author{Satoshi Koike and Laurentiu Paunescu}

\address{Department of Mathematics, Hyogo University of Teacher Education,
942-1 Shimokume, Kato, Hyogo 673-1494, Japan}
\email{koike@hyogo-u.ac.jp} 
\address{School of Mathematics and Statistics, University of Sydney, Sydney, 
NSW, 2006, Australia}
\email{laurentiu.paunescu@sydney.edu.au}


\date{\today}

\begin{document}



\begin{abstract}
In this paper we investigate the behaviour of the geometric directional bundles, associated to arbitrary subsets in $\R^n$, under bi-Lipschitz homeomorphisms, and give conditions under which their bi-Lipschitz type is preserved. The most general sets we consider satisfy the sequence selection property (SSP) and, consequently, we  investigate the behaviour of such sets under bi-Lipschitz homeomorphisms as well. In particular, we show that the  bi-Lipschitz images of a subanalytic sets generically satisfy  the (SSP) property.

\end{abstract}

\maketitle

\section{Introduction}

Given  $A$  a subset of $\R^n$
and  $p \in \overline{A}$, 
where $\overline{A}$ denotes the closure of $A$ in $\R^n$,  one may associate
 the {\em direction set or directional tangent cone} $D_p(A)$ of $A$ at $p$ (see  the Preliminary Section below). This object has been extensively studied in the literature (see for instance  \cite{Sa}, \cite{Sa1}, \cite{kp1} - \cite{kp6}), especially its behaviour under bi-Lipschitz mappings.
In this paper, for a better understanding of bi-Lipschitz geometry of sets, we employ the natural notion of geometric directional bundle at $p \in \overline{A}$ (see \cite{kp6}), thus taking into account the limit behaviour of the directional tangent cones at $q \in \overline{A}$, $q$ near $p$, and analise its behaviour under bi-Lipschitz mappings. Unexpectedly, the geometric directional bundle, is merely a $C^1$ invariant rather than a bi-Lipschitz one, hinting to a subtle geometric difference between $C^1$ diffeomorphisms and bi-Lipschitz homeomorphisms.
Despite the fact that the directional tangent cones are, in a way, bi-Lipschitz invariants, we give an example  showing that this is not the case for the geometric directional bundles (see the end of \S 4).

In $\S2$ we give the definitions of the directional tangent cone and of the sequence selection property, denoted by (SSP). We also give the definition of the geometric directional bundle and give some illustrative examples in $\S3$. In our approach to Lipschitz geometry, using directional properties,  we employ  a basic property stated in Lemma 4.1 and Remark 4.2 (also discussed  in \cite{kp1}-\cite{kp4}). In $\S4$  we discuss a similar property for the geometric directional bundle, described by the Theorem 4.7. Subsequently, we analise the (SSP) property at points in the image of a given sub-analytic set by a bi-Lipschitz homeomorphism, using the Rademacher theorem in $\S5$. The main result in this section is Theorem 5.4. 

In the final section, to better understand the Lipschitz type of a given set, we introduce several natural equivalence relations on specific collection of sets.

\bigskip
\section{preliminaries}
\subsection{Direction set}

We first recall the definition of the direction set.

\begin{defn}\label{direction}
Let $A$ be a subset of $\R^n$, 
and let $p \in \R^n$ such that $p \in \overline{A}$, 
where $\overline{A}$ denotes the closure of $A$ in $\R^n$.
We define the {\em direction set} $D_p(A)$ of $A$ at $p$ by
$$
D_p(A) := \{a \in S^{n-1} \ | \
\exists  \{ x_i \} \subset A \setminus \{ p \} ,
\ x_i \to p \in \R^n  \ \text{s.t.} \
\frac{x_i - p}{\| x_i - p \| } \to a, \ i \to \infty \}.
$$
Here $S^{n-1}$ denotes the unit sphere centred at $0 \in \R^n$.

For a subset $D_p(A) \subset S^{n-1}$, we denote by $LD_p(A)$
the semi-cone of $D_p(A)$ with $0 \in \R^n$ as the vertex, and call it the {\em real tangent cone} of $A$ at $p$:
$$
LD_p(A) := \{ t a \in \R^n\ | \ a \in D_p(A), \ t \ge 0 \}.
$$
For $p \in \overline{A}$ we put $L_pD(A) :=p+ LD_p(A)$, and call it the {\em geometric tangent cone} of $A$ 
at $p$.
In the case where $p = 0 \in \R^n$, we write
$D(A) := D_0(A)$ and $LD(A) := L_0D(A)$ for short.
\end{defn}

Concerning the direction set, we introduced  in \cite{kp1} 
the notion of the sequence selection property, denoted by $(SSP)$ for short.
It is a very useful property to study the directional properties of singular spaces. 
For instance, see \cite{kp3, kp4, kp5}. 
Let us recall the notion of $(SSP)$.

\begin{defn}
Let $A \subset \R^n$ and  $p \in \R^n$  such that $p \in \overline{A}$. 
We say that $A$ satisfies {\em condition} $(SSP)$ {\em at} $p$, if for any sequence of points $\{ a_m\}$
of $\R^n$ tending to $p \in \R^n$ such that $\lim_{m \to \infty} \frac{a_m - p}{\| a_m - p\|} \in D_p (A)$,  
there is a sequence of points $\{ b_m\} \subset A$ such that 
$$
\| a_m - b_m\| \ll \| a_m - p\|, \ \| b_m - p\| .
$$
\end{defn}

As mentioned in the Introduction of \cite{kp3}, $(SSP)$ is shared by a large class of sets,
in particular subanalytic sets and cones. 

\subsection{Geometric directional bundle}

From an attempt to produce global Lipschitz invariants, we introduced the notion of 
geometric directional bundle in \cite{kp6}. 
Let us recall the notion here.

\begin{defn}\label{gdirectionset} 
Let $A \subset \R^n$, and let $W \subset \R^n$ such that 
$\emptyset \ne W \subset \overline{A}$. 
We define the {\em direction set} $D_W(A)$ of $A$ over $W$ by
$$
D_W (A) := \bigcup_{p \in W} (p, D_p(A)) \subseteq W \times S^{n-1}\subseteq \R^n\times S^{n-1}.
$$

We call 
$$\mathcal{GD}_W(A) := W\times S^{n-1}\cap \overline{D_A(A)},$$
the {\em geometric directional bundle} of $A$ over $W,$
where $\overline{D_A(A)}$ denotes the closure  of $D_A(A)$
in $\R^n \times S^{n-1}$.

Let $\Pi : \R^n \times S^{n-1} \to S^{n-1}$ be the canonical
projection defined by $\Pi(x,a)=a$.
In the case where $W = \{ p \}$ for $p \in \overline{A}$, we write
$$\mathcal{GD}_p(A) := \Pi(\mathcal{GD}_{\{p\}}(A))\subseteq S^{n-1}.$$  
This is the set of all possible limits of directional sets $D_q(A), q\in A, q\to p.$

We consider the half cone of $\mathcal{GD}_p (A)$ with $0\in \R^n$ as the vertex
$$
L\mathcal{GD}_p(A) :=  \{ t v \in \R^n \ | \ v \in \mathcal{GD}_p(A), \ t \ge 0 \}.
$$ We call it the {\em real tangent bundle cone} of $A$ at $p$.

Similarly we define the {\em geometric bundle cone} of $A$ at $p,$ to be its translation by $p$
$$
L_p\mathcal{GD}(A) := p+ L\mathcal{GD}_p(A) .
$$ 
In the case where $p = 0 \in \R^n$, we simply write
$\mathcal{GD}(A) := \mathcal{GD}_0(A)$, and 
$L\mathcal{GD}(A)  := L_0\mathcal{GD}(A)$.
\end{defn}

In \cite{kp7} we discussed the stabilisation of the geometric directional bundle 
as an operation acting on singular spaces.

\subsection{Local dimension of a subanalytic set}

Let $A \subset \R^n$ be a subanalytic set and let $p \in A$. 
It is known by H. Hironaka \cite{hironaka} that a subanalytic set admits a locally finite 
Whitney stratification with analytic strata. 
Therefore there exist a positive number $\epsilon_0 > 0$ and an integer $k$ 
with $0 \le k \le n$ such that for any $0 < \epsilon < \epsilon_0$, 
$A \cap B_{\epsilon}(p)$ is a $k$-dimensional subanalytic set, 
where $B_{\epsilon}(p)$ denotes an $\epsilon$-neighbourhood of $p$ in $\R^n$. 
(For the Whitney stratification, see H. Whitney \cite{whitney1, whitney2}.) 
We call this k {\em the local dimension of $A$ at $p$}.

Using the local dimension, we introduced a dimensional condition in \cite{kp7}. 
We recall it here.

\begin{defn}
Let $A \subset \R^n$ be a $k$-dimensional subanalytic set ($0 \le k \le n$). 
We say that $A$  is {\em genuinely $k$-dimensional}, 
if for any $p \in A$, the local dimension of $A$ at $p$ is $k$. 
\end{defn}

\begin{example}\label{Whitneyumbrella}[not genuinely 2-dimensional]

Let $A \subset \R^3$ be the Whitney umbrella defined by 
$$
A := \{ (x,y,z) \in \R^3 \ | \ x^2 = z y^2 \}. 
$$
This is a $2$-dimensional algebraic set.
Let $A_1$ (respectively $A_2$) be the set of points of $A$ 
at which the local dimension is $1$ (respectively $2$). 
Then we have 
$$
A_1 = A \cap \{ (x,y,z) \in \R^3 \ | \ z < 0 \}, \ 
A_2 = A \cap \{ (x,y,z) \in \R^3 \ | \ z \ge 0 \}. 
$$
Therefore the Whitney umbrella $A$ is not genuinely $2$-dimensional.
We note that $\overline{A_1} \cap A_2 = \{ (0,0,0)\}$ and 
$A_1 \cap \overline{A_2} = \emptyset$. 
\end{example}

\subsection{The Rademacher theorem}

In this subsection we recall the Rademacher theorem about the differentiability 
of Lipschitz maps.

\begin{thm}\label{rademacher} (Rademacher's theorem) 
Let $U$ be an open subset of $\R^n$ and let $f : U \to \R^m$ be a Lipschitz map.
Then $f$ is differentiable almost everywhere in $U$, that is, the points at which 
$f$ is not differentiable form a set of Lebesque measure zero.
\end{thm}

For a proof of the above theorem and some useful corollaries of it, 
 see for instance  chapter 3 in L.C. Evans and R.F. Gariepy \cite{rade}.

\bigskip
\section{Examples of geometric bundle cones}

In this section we determine the geometric tangent cones and the geometric bundle cones
of several  examples.

\vspace{3mm}

\noindent {\bf (1)} Let
$$
A_1 := \{ (x,y,z) \in \R^3 \ | \ x^2 + y^2 = z^2 \}, \ A_2 := \{ (x,y,z) \in \R^3 \ | \ x^2 + y^2 = z^4 \}.
$$
Concerning the geometric tangent cones, we have $LD(A_1) = A_1$ and $LD(A_2) = \{ z-{\text axis}\}$. 
Therefore we have  $\dim LD(A_1) = 2$ and $\dim LD(A_2) = 1$. 
It follows from Theorem 1.1 in \cite{kp1} that $A_1$ and $A_2$ are not bi-Lipschitz 
equivalent as germs at $0 \in \R^3$.  

Concerning the geometric bundle cones, we have
$$
L\mathcal{GD}(A_1) =  \{ (x,y,z) \in \R^3 \ | \ z^2 \le x^2 + y^2  \}, \ \ L\mathcal{GD}(A_2) = \R^3 
$$
(see Example 2.4 (1) in \cite{kp7} for the details).

\vspace{3mm}

\noindent{\bf (2)}[dimension invariance]

Let $P_1 = (1,1,1)$, $P_2 = (1,-1,1)$, $P_3 = (-1,-1,1)$ and $P_4 = (-1,1,1)$ be points 
on the plane $\{ z = 1\}$ in $\R^3$, and let $Q \subset \R^3$ be the square 
consisting of segments $P_1 P_2$, $P_2 P_3$, $P_3 P_4$ and $P_4 P_1$. 
We denote by $B$ the real cone with vertex $(0,0,0)$ and base $Q$, namely 
$$
B := \{ t a \in \R^3\ | \ a \in Q, \ t \ge 0 \}.
$$
If we put  
$$
A_1^+ := \{ (x,y,z) \in \R^3 \ | \ x^2 + y^2 = z^2 \ \& \ z \ge 0 \},
$$
then we can see that $B$ is bi-Lipschitz equivalent to $A_1^+$.

Concerning the geometric tangent cone, we have $LD(A_1^+) = A_1^+$ and $LD(B) = B$.
Therefore $\dim LD(A_1^+) = \dim LD(B) = 2$. 
We have shown a general result on the dimension of the geometric tangent cone 
(see \cite{kp1} Theorem 1.1).

\begin{thm}\label{lipinvariant} 
Let $h :(\R^n,0) \to (\R^n,0)$ be a bi-Lipschitz homeomorphism.
If $A$, $h(A)$ are subanalytic set-germs at $0 \in \R^n$, 
then $\dim LD(h(A)) = \dim LD(A)$.
\end{thm}

Let $\tau_1$, $\tau_2$, $\tau_3$ and $\tau_4$ be hyperplanes in $\R^3$  defined by
$$
x - z = 0, \ \ y + z = 0, \ \ x + z = 0, \ \ y - z = 0,
$$
respctively.
Then we have 
$$
L\mathcal{GD}(A_1^+) =  \{ (x,y,z) \in \R^3 \ | \ z^2 \le x^2 + y^2  \}, \ \ 
L\mathcal{GD}(B) = \tau_1 \cup \tau_2 \cup \tau_3 \cup \tau_4 .
$$
It follows that $\dim L\mathcal{GD}(A_1^+) = 3$ but $\dim L\mathcal{GD}(B) = 2$. 
We note that $L\mathcal{GD}(A_1^+) $ and $L\mathcal{GD}(B)$ are semialgebraic sets in $\R^3$.  
Therefore, even in the semialgebraic case, $\dim L\mathcal{GD}(A)$ is not a bi-Lipschitz invariant.

\vspace{3mm}

\noindent {\bf (3)} Let $f : (\R^n,0) \to (\R,0)$ be a continuous function defined in a neighbourhood of $0 \in \R^n$. 
We say that $f$ satisfies the {\em flat condition} at $0 \in \R^n$ if there exit $e > 0$, $C > 0$ and $\delta > 0$ such that 
$|f(x)| \le C \| x\|^{1 + e}$ for $\| x\| < \delta$. 
We denote by $A$ the graph of $f$.
Then we have the following claims.

\vspace{3mm}

\noindent {\em Claim 1.} If a continuous function $f : (\R^n,0) \to (\R,0)$ satisfies the flat condition at $0 \in \R^n$,
then we have $LD(A) = \R^n \times \{ 0\} \subset \R^{n+1}$.

\vspace{3mm}

\noindent {\em Claim 2.} If a $C^1$ function $f : (\R^n,0) \to (\R,0)$ satisfies the flat condition at $0 \in \R^n$,
then we have $L\mathcal{GD}(A) = \R^n \times \{ 0\} \subset \R^{n+1}$.

\vspace{3mm}

The following example shows that we cannot drop the assumption of $C^1$ in Claim 2.

\begin{example}
Let us recall the real Brian\c{c}on-Speder family (\cite{BrSp}) $f_t : (\R^3,0) \to (\R,0)$, 
$t \in J = ( - (1 + \epsilon ), 1 + \epsilon )$, defined by 
$$
f_t (x,y,z) := z^5 + tzy^6 + y^7 x + x^{15}, 
$$
where $\epsilon$ is a sufficiently small positive number. 
Therefore $\{ f_t \}$ is a family of polynomial functions with isolated singularities. 

Let $A := f_0^{-1}(0)$. 
Then we can see that $A$ is the graph of a semialgebraic, continuous function $g : (\R^2,0) \to (\R,0)$ 
defined by $g(x,y) = - (y^7x + x^{15})^{\frac{1}{5}}$, and it satisfies the flat condition at $0 \in \R^2$. 
As mentioned in \cite{Ko}, $g$ is differentiable at $0 \in \R^2$ but not $C^1$. 
By Claim 1, we have $LD(A) = \R \times \R \times \{ 0\} \ (= \R^2 \times \{ 0\})$. 
The gradient of $f_0$ is as follows:
$$
\grad f_0 (x,y,z) = (y^7 + 15 x^{14}, 7y^6x, 5z^4).
$$

Let $\ell := \{ (x,y,z) \in \R^3 \ | \ x = z = 0 \}$, namely the $y$-axis in $\R^3$. 
Then we have $\ell \subset A$ and $\grad f_0 (0,y,0) = (y^7,0,0)$. 
Therefore we have 
$$
\lim_{y \to 0} \frac{\grad f_0 (0,y,0)}{\| \grad f_0 (0,y,0) \| } = (1,0,0) .
$$
It follows that 
\begin{eqnarray}
\{ 0 \} \times \R \times \R \subset L\mathcal{GD}(A). 
\end{eqnarray}

Let {\em m} $:= \{ (x,y,z) \in \R^3 \ | \ y = - x^2 \ \& \ z = 0 \}$. 
Then we have {\em m} $\subset A$ and $\grad f_0 (x,- x^2,0) = (14x^{14},7x^{13},0)$. 
Therefore we have 
$$
\lim_{x \to 0} \frac{\grad f_0 (x,-x^2,0)}{\| \grad f_0 (x,-x^2,0) \| } = (0,1,0) .
$$
It follows that 
\begin{eqnarray}
\R \times \{ 0 \} \times \R \subset L\mathcal{GD}(A). 
\end{eqnarray}

Let us denote by $S$ the $2$-dimensional unit sphere centred at $0 \in \R^3$. 
Then $S \cap A$ is homeomorphic to a circle. 
Note that if $(x_0 , y_0, z_0 ) \in S \cap (A \setminus (\ell \ \cup$ {\em m})), 
then we have $z_0 \ne 0$. 
For any $(x_0 , y_0, z_0 ) \in S \cap (A \setminus (\ell \ \cup$ {\em m})), 
we define a curve $\lambda$ by $\lambda (s) = (x_0 s, y_0 s^2 , z_0 s^3)$, $0 < s \le 1$. 
Then we have $\lambda \subset A$ and 
$$
\grad f_0 (x_0 s, y_0 s^2, z_0 s^3) = ((y_0^7 + 15 x_0^{14})s^{14}, 7y_0^6 x_0 s^{13}, 5z_0^4 s^{12}).
$$
Since $z_0 \ne 0$, we have
$$
\lim_{s \to 0} \frac{\grad f_0 (x_0 s, y_0 s^2, z_0 s^3)}{\| \grad f_0 (x_0 s, y_0 s^2, z_0 s^3)\|} 
= (0, 0, 1).
$$
It follows that 
\begin{eqnarray}
\R \times \R \times \{ 0 \} \subset L\mathcal{GD}(A). 
\end{eqnarray}

By (3.1), (3.2) and (3.3), we see that 
$$
\R \times \R \times \{ 0 \} \varsubsetneqq L\mathcal{GD}(A). 
$$
\end{example}

\bigskip

\section{Geometric directional bundle properties related to (SSP)}

One of the important directional properties related to $(SSP)$ is the following lemma. 

\begin{lem}\label{sspcone} $($\cite{kp5}, Lemma 3.4$)$ 
Let $A \subset \R^n$ and  $p \in \R^n$  such that $p \in \overline{A}$, 
and let $h : (\R^n ,p) \to (\R^n ,h(p))$ be a 
bi-Lipschitz homeomorphism. 
Suppose that $A$ satisfies condition $(SSP)$ at $p$. 
Then for some 
bi-Lipschitz homeomorphism 
$\tilde h$  we have 
$$
\tilde h(LD_p (A)) \subseteq LD_{h(p)} (h(A)).
$$
\end{lem}

\begin{rem}\label{remark11} 
In the above lemma we have equality if $h(A)$ 
satisfies  the assumption of $(SSP)$ at $h(p)$, see Theorem 3.6 in \cite{kp5}.
\end{rem} 

On the other hand, the lemma in the case of geometric directional bundle  
corresponding to Lemma \ref{sspcone} has the following formulation.

\begin{lem}\label{ssphatlemma} 
Let  $A \subset \R^n, 0\in \overline A$  
and let $h : (\R^n ,0) \to (\R^n ,0)$ be a bi-Lipschitz homeomorphism. 
Let $p_m\in \overline A\to 0 $ and suppose that $A$ satisfies condition $(SSP)$ at all $p_m\in \overline A$. Then  for some bi-Lipschitz homeomorphism $\tilde h$  we have 
$$
\tilde h(\lim LD_{p_m}(A)) \subseteq \lim LD_{h(p_m)}h(A).
$$

The  inclusion is an equality if $h(A)$ satisfies condition $(SSP)$ at $h(p_m ) \in h(A),$ for any  $m \in \N.$ 
Moreover if $p_m\in A\setminus \{0\}$ then we have
$$\tilde h(\lim LD_{p_m}(A)) \subseteq \lim LD_{h(p_m)}h(A) \subseteq L\mathcal{GD} (h(A)).$$
\end{lem}

\begin{rem}\label{remark12} 
The above lemma is a generalisation of the Lemma \ref{sspcone} above (take $p_m=0,$  for any $m\in \N$).
\end{rem} 

\begin{proof}
The proof of the above result is similar to the proof of the Lemma \ref{sspcone} above. Indeed for each $p_m\in \overline A,$ via Arzela-Ascoli (see \cite{Sa}) there is a bi-Lipschitz homeomorphism $d_{p_m} h_m$, with the same Lipschitz constant as $h$, such that $ d_{p_m}h_m(LD_{p_m}(A)) \subseteq  LD_{h(p_m)}h(A),$	with equality if $h(A)$ satisfies condition $(SSP)$ at $h(p_m)$. Now again by using the Arzela-Ascoli result we can claim that (modulo a subsequence) $d_{p_m}h_m$ converges to some bi-Lipschitz homeomorphism $\tilde h$ which will have the required property. Indeed it is not hard to see that if $v_m\in LD_{p_m}(A), v_m\to v\in \lim LD_{p_m}(A)$ then $\tilde h(v)=\lim d_{p_m} h_m(v_m)\in \lim LD_{h(p_m)}h(A) $.
\end{proof}

Let $(\overline A, 0)\subset (\R^n,0)$ and assume that there is a family of subspaces $\pi_i\subset \R^n, i\in I$ such that $\overline {\cup_{i\in I} \pi_i }=L\mathcal{GD} (A)$. Note that in the case when $A$ is subanalytic, $LDA$ is a subset of the union of the limits of subspaces 
$\lim LD_pA, \,\, p\in A, p\to 0$ i.e. a subset of $L\mathcal{GD} (A)$.
For each $\pi_i, i\in I$ we assume that there is  a sequence of points $p_m(i)\in A, p_m(i)\to 0, m\to\infty$ such that 
there is a bi-Lipschiz mapping 
$h_{m,i}$, such that $ h_{m,i}( LD_{p_m(i)}(A))\subseteq  LD_{h_{m,i}(p_m(i))}h(A).$ 
Moreover we assume that  $h_{m,i}\to h_i, m\to \infty$,  satisfying the following 
$ LD_{p_m(i)}(A)\to \pi_i, m\to \infty, \, \, h_i(\lim_m  LD_{p_m(i)}(A))=
h_i(\pi_i)\subseteq \lim_m LD_{h_{m,i}(p_m(i))}h(A) \subseteq L\mathcal{GD} (h(A)).$

If for each $i, j\in I$ we have $h_i=h_j=\tilde h$ then we can conclude that 
$$
\tilde h (\overline {\cup_{i\in I} \pi_i }) =\tilde h (L\mathcal{GD} (A)) \subseteq L\mathcal{GD} (h(A)).
$$

\begin{rem}
If $A$ satisfies condition $(SSP)$ at all $p\in \overline A$ in a small enough neighbourhood of  the origin, by using  Lemma \ref{ssphatlemma} one can create the above $h_{m,i}\in I$ as Arzela-Ascoli  derivatives of the initial $h$, and denote them by $d_{p_m(i)}h$, but we cannot always assume that $h_i$ (also obtained by the Arzela-Ascoli process, given by the proof of Lemma \ref{ssphatlemma}) are independent on $i\in I$, thus we need to consider special classes of bi-Lipschitz homeomorphisms to compare the directional bundles.
\end{rem}

\begin{defn}\label{specialdir}
Let $h : (\R^n ,0) \to (\R^n ,0)$ be a bi-Lipschitz homeomorphism and let  $A \subset \R^n, 0\in \overline A.$  We say that $h$ is  $A-$\emph{directional} (at $0\in \overline A,$) whenever there is a  subset $D$ of $A$
such that $0\in \overline {D}$ and for $p \in D,$
$$
d_{p}h\to \tilde h \ \text{as} \ p\to 0, \ d_ph(LD_p (A)) \subseteq LD_{h(p)} (h(A)) \ 
\text{and} \ L\mathcal{GD} (A)\subseteq  \overline {\cup_{p \in U\cap D}LD_pA},
$$ for all neighbourhoods $U$ of $0$.
\end{defn}

For example, the bi-Lipschitz mapping $h: \R^2 \to \R^2, \ h(x,y)=(x,y+|x|)$, is $A-$\emph{directional}  but $h^{-1}$ is not $h^{-1}(A)-$\emph{directional}  for $A=(-\infty,\infty).$ Obviously every $C^1$ diffeomorphism is $A-$\emph{directional}  for any subset $A$.

In Definition \ref{specialdir} one may ask the convergence condition to be independent on $p$ only for the restrictions to $A$, 
that is $d_{p}(h/A)\to \tilde h, p\to 0.$ It follows that for such $\tilde h$ we have 
$$\tilde h (L\mathcal{GD} (A)) \subseteq L\mathcal{GD} (h(A)).$$

\begin{thm}
Let  $A \subset \R^n, 0\in \overline A$ and  let $h : (\R^n ,0) \to (\R^n ,0).$  
Assume that $A$ satisfies condition $(SSP)$ at all $p\in \overline A$ in a small enough neighbourhood of the origin and $L\mathcal{GD}_{0} (A)\subseteq  \overline {\cup_{p\in U\cap D} LD_pA}$, $U$ arbitrary neighbourhood of the origin.
Then 
$h$ is an $A-$directional bi-Lipschitz homeomorphism provided that there is a subset $D$ of $A$, $0\in \overline {D},$ such that 
$\lim_{p \to 0} d_{p}h = \tilde h$ is independent of $p\in D$. 
 In this case we have 
$$\tilde h (L\mathcal{GD} (A)) \subseteq L\mathcal{GD} (h(A)).$$
Moreover if $h(A)$ also satisfies condition $ (SSP)$ at every point near the origin and also $h^{-1}$ is $h^{-1}(A)-$directional, then the above inclusion becomes equality.
\end{thm}

%
%

\begin{rem}
In the case $\tilde h$ above is also resulting from the Arzela-Ascoli process i.e.  $\tilde h=d_0h$ then we know that we also have Lemma \ref{sspcone} 
$$
\tilde h(LD (A)) \subseteq LD (h(A)).
$$ or equality if  $h(A)$ also satisfies condition $ (SSP)$ at the origin. 

\end{rem}
The following example shows that the geometric directional bundle is not a bi-Lipschitz invariant.
\begin{example} Let  $h:\R^2\to\R^2$ defined by the following expression 
$h(x,y)=(2x+(x^2+y^2)\sin\frac{1}{\sqrt{x^2+y^2}},y)$, and $A:=\{(0,t)|t\in \R\}$. This $h$ is analytic except the origin, differentiable everywhere and bi-Lipschitz, but $h(0,t)=(t^2\sin\frac{1}{|t|},t)$ and $A$ and $h(A)$ have different geometric directional bundles. This shows that the notion of geometric directional bundle is not a bi-Lipschitz invariant, but rather a strict $C^1$ invariant. Note that both $A$ and $h(A)$ satisfy the (SSP) property.
\end{example}

\bigskip
\section{(SSP) of the image by a bi-Lipschitz homeomorphism}

In the sequel we will use the following  terminology. 
Let $\{ e_i\}$ be a sequence of points of $\R^n \setminus \{ 0\}$ tending to $0 \in \R^n$ 
such that $\frac{e_i}{\| e_i\|}$ tends to $e \in D\R^n$. 
Then we say that $e_i$ {\em induces the direction} $e$ or 
$e$ is {\em the direction given by} $e_i$.

\begin{lem}
Let $h:(\R^n,0)\to (\R^n,0)$ be a bi-Lipschitz homeomorphism and $v\in S^{n-1}$ such that $\frac{h(tv)}{t}$ is convergent as $t\to 0, t>0.$ 
Then, if we put $ l_v:=\{tv | t>0\}$ it follows that $Dh(l_v)=\{c\}$ is a point.
Moreover if  $c_i=h(b_i)\to 0, i\in \N$ induces the direction $c$,  then $b_i$ induces the direction $v$. In particular $h(l_v)$ satisfies (SSP) as well.
\end{lem}
\begin{proof} Using our assumption it is easy to see that $h(t_iv)$ induces a unique direction $c$, that is $Dh(l_v)=\{c\}$ .
Let $\{ c_i=h(b_i), i \in \N\}, b_i\to 0 \in \R^n$ be an arbitrary sequence of points of $\R^n$ tending to $0 \in \R^n$ such that $\frac{c_i}{\|c_i\|}$ tends to the direction $\{c \}=Dh(l_v)$ ($c_i$ induces the direction $c$).

By definition, $c$ is also given by $\frac{h(a_i)}{\|h(a_i)\|}$ for some sequence $\{ a_i=t_iv\}$ of $l_v$ tending to $0  \in \R^n$.
Let $dh (w)$ be the bi-Lipschitz mapping  given by Arzela-Ascoli (we call it an Arzela-Ascoli differential) applied to the sequence 
$h_i(w):= \frac{h(\| b_i\|w)}{\| b_i\|}$. 
By assumption, $c$ is the direction given by $\frac{h(a_i)}{\|h(a_i)\|}$ and also the direction induced by $\frac{h(\|a_i \| v)}{\| a_i\|}.$ 

Indeed one can write 
$$
\frac{h(a_i)}{\|h(a_i)\|} = \frac{h(\|a_i\| \frac{a_i}{\|a_i\|})}{\|h(a_i)\|} = (\frac{h(\|a_i\| \frac{a_i}{\|a_i\|})}{\|a_i\|}) \times (\frac{\|a_i\|}{\|h(a_i)\|}),
$$ 
and this induces the same direction as $\frac{h(\|a_i \| v)}{\| a_i\|}$.

By assumption, the limit $\frac{h(tv)}{t}$ is independent of $t$ so we can use $\| b_i \|$, namely $\frac{h(\| a_i\| v)}{\|a_i \|} \to dh(v)$. 
Accordingly the direction $c$ is both $dh(\beta b), \, b$ a direction induced by $b_i$,  and $dh(\alpha v)$, $\alpha, \beta >0$ given by the limits of 
$\frac{\|a_i\|}{\|h(a_i)\|}$ and $\frac{\|b_i\|}{\|h(b_i)\|}$ respectively. 
Because $dh$ is bi-Lipschitz itself it will follow that $\{b\}= Dl_v=\{v\}$.
\end{proof}

\begin{prop}\label{imSSP}
Let $h:(\R^n,0)\to (\R^n,0)$ be a bi-Lipschitz homeomorphism such that for any $v\in DA, h(tv)/t\to dh(v)$ as $t\to 0.$
If $A$ satisfies condition $(SSP)$ then $h(A)$ satisfies condition $(SSP)$ and we have  $LDh(LDA)=LDh(A)=dh(LDA)$. 
Here $dh$ can be any Arzela-Ascoli differential $dh$ of $h$, in particular it will be bi-Lipschitz itself.
\end{prop}

\begin{proof} 
Let $\{ c_i=h(b_i), i \in \N\}, b_i\to 0 \in \R^n$ be an arbitrary sequence of points of $\R^n$ tending to $0 \in \R^n$ such that $\frac{c_i}{\|c_i\|}$ tends to a direction 
$c\in Dh(A)$ ($c_i$ induces the direction $c$). Every direction $c\in Dh(A)$ can be realised as a direction coming from $l_v, v\in DA$ conveniently chosen. 
Then $b$ the direction induced by $b_i$ is $v\in DA$ by the above Lemma.

Since $A$ satisfies condition $(SSP)$, it follows that there are $d_i\in A$ such that $\|b_i-d_i\| \ll \|b_i\|$ which in turn it will show that 
$\| c_i -h(d_i)\| \ll \| c_i \|$. 
Hence $h(A)$ satisfies condition $(SSP)$ itself.
The rest is a consequence of the Lemma \ref{sspcone}.

\end{proof}

Using the notion of a genuinely $k$-dimensional subanalytic set defined in \S 2.4,
as a corollary of Proposition \ref{imSSP}, we have the following lemma.

\begin{lem}\label{denselem}
Let $A \subset \R^n$ such that $0\in \overline A,$ and let $h : (\R^n ,0) \to (\R^n ,0)$ be a bi-Lipschitz homeomorphism.
If $A$ is a genuinely $i$-dimensional subanalytic set, $0 \le i \le n$, then $h(A)$ satisfies condition $(SSP)$ on a dense subset.
\end{lem}

\begin{proof}
If $i = 0$, then $A = h(A) = \{ 0 \}$ as germs at $0 \in \R^n$. 
Therefore $h(A)$ satisfies condition $(SSP)$ on $h(A)$. 
After this, we assume that $1 \le i \le n$.

Let us denote by $R_A$ the regular set of $A$. 
Since $A$ is a genuinely $i$-dimensional subanalytic set, $R_A$ is an open dense subset of $A$. 
On the other hand, the Rademacher theorem applied to the regular set $R_A$ gives 
a dense set $S$ in $R_A$ where $h|_{R_A}$ is differentiable. 
Note that $S$ is dense also in $A$ since $A$ is genuinely $i$-dimensional subanalytic.
At the points of $S$ any Arzela-Ascoli derivative will have the property from the above proposition. 
This implies that if $A$ is a genuinely $i$-dimensional subanalytic set, then $h(A)$ will satisfy condition $(SSP)$ on a dense subset $h(S)$.
\end{proof}

We can generalise the above lemma as follows.

\begin{thm}\label{density}
Let $A \subset \R^n$ such that $0\in \overline A,$ and let $h : (\R^n ,0) \to (\R^n ,0)$ be a bi-Lipschitz homeomorphism.
If $A$ is subanalytic, then $h(A)$ satisfies condition $(SSP)$ on a dense subset.
\end{thm}

\begin{proof}

Let $A \subset \R^n$ be a $k$-dimensional subanalytic set, $0 \le k \le n$, 
such that $0 \in \overline{A}$. 
Let us recall the decomposition of $A$ into genuinely $i$-dimensional subanalytic subsets ($0 \le i \le k$) 
introduced in \cite{kp7}.
We define the following subsets of $A$,
$$
A_i := \{ p \in A : \text{the local dimension of} \ A \ 
\text{at} \ p \ \text{is} \ i \}
$$
for $0 \le i \le k$. 
Locally at $0 \in \R^n$ $A$ admits a finite Whitney stratification $\mathcal{S}(A)$. 
We are assuming the frontier condition (in the sense of J.N. Mather \cite{mather}) for this stratification.
Therefore the local dimension is constant on each stratum of $\mathcal{S}(A)$.
It follows that if $A_i$, $0 \le i \le k$, is non-empty, $A_i$ is a genuinely  $i$-dimensional 
subanalytic subset of $A$.
In particular, since $A$ is a $k$-dimensional subanalytic set, $A_k$ is non-empty and 
a genuinely $k$-dimensional subanalytic set.

Set
$$
\Lambda_0 (A) := \{ \ 0 \le i \le k \ | \ 0 \in \overline{A_i} \} .
$$
Then we can locally express $A$ as a disjoint union  around $0 \in \R^n$ as follows:
\begin{equation*}\label{localunion}
A \cap B_{\epsilon}(0) 
= \bigcup_{i \in \Lambda_0 (A)} A_i \cap B_{\epsilon}(0) 
\end{equation*}
for a sufficiently small $\epsilon > 0$. 

We denote by $R_{A_i}$ the regular set of $A_i$ for $i \in \Lambda_0 (A)$ 
like $R_A$ in the proof of the above lemma. 
Let us apply Lemma \ref{denselem} to each $A_i$.
Then $h(A_i )$ satisfies condition $(SSP)$ on a dense subset $h(S_i)$ 
for $i \in \Lambda_0 (A)$. 
Here each $S_i$ is given like $S$ in the above lemma.
Therefore we have 
\begin{equation}\label{subset}
h(S_i) \subset h(R_{A_i}) \ \text{for} \ i \in \Lambda_0 (A).
\end{equation}

Let $i, j \in \Lambda_0 (A)$ such that $i < j$.
Even if $\overline{h(A_i)} \cap h(A_j) \ne \emptyset$, we have $(\overline{h(A_i)} \cap h(A_j)) \cap h(S_j ) = \emptyset$ 
from (\ref{subset}). 
This means that $h(A)$ satisfies condition $(SSP)$ at any $p \in h(S_j)$ for $j \in \Lambda_0 (A)$. 
Set $S = \bigcup_{i \in \Lambda_0 (A)} S_i$. 
Then, by construction, $h(S)$ is a dense subset of $h(A)$ and $h(A)$ satisfies condition $(SSP)$ on $h(S)$.
This completes the proof.
\end{proof}

\begin{rem}
The images of subanalytic sets $A$ of bi-Lipschitz homeomorphisms $h$ will have for a dense subset $B\subseteq h(A)$ the same directional cone and the same geometric directional bundle i.e. $LD_qh(A)=L\mathcal{GD}_{q} (h(A)), q\in B$ (as this property clearly holds for subanalytic sets). Note that in general there is no obvious relation between $LD_q(A) \,\,\text{and} \,\,L\mathcal{GD}_{q}(A),$
nevertheless in the subanalytic case we always have $LD_q(A)\subseteq L\mathcal{GD}_{q}(A)$. In fact our results can be proved for $(SSP)$ sets $A$ which locally are finite union of bi-Lipschitz images of open subsets in Euclidean spaces.
We can also prove our results just for Lipschitz mappings with $\tilde h$ Lipschitz and mention that if $h$ is bi-Lipschitz then $\tilde h$ is also bi-Lipschitz.
\end{rem}

By Theorem 5.11 in \cite{kp2}, we have the following lemma.

\begin{lem}\label{h(LDA)}
Let $A \subset \R^n$ such that $0\in \overline A,$ and let $h : (\R^n ,0) \to (\R^n ,0)$ be a bi-Lipschitz homeomorphism.
Assume that $A$ satisfies condition $(SSP)$ at $p \in \overline{A}$.
Then $h(A)$ satisfies condition $(SSP)$ at $h(p) \in \overline{h(A)} = h(\overline{A} )$ if and only if 
$h(LD_p(A))$ satisfies condition $(SSP)$ at $h(p) \in \overline{h(A)}$.
\end{lem}

Let $A \subset \R^n$ be a subanalytic subset such that $0\in \overline A,$ 
and let $h : (\R^n ,0) \to (\R^n ,0)$ be a bi-Lipschitz homeomorphism. 
Let $S$ be a subset of $A$ given in the proof of Theorem \ref{density}.
Then, as seen in the proof, $h(S)$ is a dense subset of $h(A)$.
From Theorem \ref{density} and Lemma \ref{h(LDA)} we have the following corollary.

\begin{cor}
Let $A \subset \R^n$ such that $0\in \overline A,$ and let $h : (\R^n ,0) \to (\R^n ,0)$ be a bi-Lipschitz homeomorphism.
If $A$ is subanalytic, then $h(LD_{h^{-1}(q)}(A))$ satisfies condition $(SSP)$ 
at any $q \in h(S)$.
\end{cor}

At the end of this section we mention the meaning of differentiable on a image of a bi-Lipschitz mapping.

\begin{prop}
Let $h:A\to \R^m, A\subseteq \R^n$, $A$ open, be a bi-Lipschitz mapping, and let $\tau :\R^m\to \R^k$ be a mapping such that $\tau/h(A)$ is Lipschitz. 
Then $\tau/h(A)$ is differentiable on $h(A)$ almost everywhere on a set of positive $n-$Lebesque measure.
\end{prop}
\begin{proof} At the points where $\tau$ itself is differentiable, there is nothing to prove (note that in general those points may be disjoint of $h(A)$ if $n<m$).

From Rademacher's theorem we know that both $h$ and $\tau\circ h$ are differentiable   almost everywhere on a set of positive $n-$Lebesque measure  $B\subseteq A$. That is, for any $x_0\in B$  there are linear mappings $L:\R^n\to \R^m, M:\R^n\to \R^k$ such that
$$h(x)-h(x_0)=L(x-x_0)+o_1(x-x_0), \, \tau\circ h(x)-\tau\circ h(x_0)=M(x-x_0)+o_2(x-x_0).$$

Because $h$ is bi-Lipschitz it follows that $L$ is injective.
Indeed $$1\sim \frac{h(x)-h(x_0)}{\|x-x_0\|}\sim L(\frac{x-x_0}{\|x-x_0\|})$$ so $L$ must have trivial kernel. Take $N:\R^m\to \R^n $ a linear left inverse of $L$.
Then we can see that $$N(h(x)-h(x_0))=N(L(x-x_0)+o_1(x-x_0))=x-x_0+o_3(x-x_0).$$
Accordingly $$\tau\circ h(x)-\tau\circ h(x_0)=M(x-x_0)+o_2(x-x_0)=M(N(h(x)-h(x_0))-o_3(x-x_0))$$
which gives $$\tau\circ h(x)-\tau\circ h(x_0)=MN(h(x)-h(x_0))+o_4(x-x_0)=MN(h(x)-h(x_0))+o_5(h(x)-h(x_0)),$$ which proves the claim on $h(B)$.
\end{proof}

\begin{rem} 
There is a Lipschitz homeomorphism with the property that its derivative vanishes on a dense set.
Let $f: [0,1]\to \R, $ be given by the absolute convergent series $ f(x) = \sum_{n \geq 1}\frac{(x-a_n)^{1/3}}{2^n}$ 
($\{a_n\in [0,1] |n\in \N\}$  dense in $[0,1]$, Pompeiu's homeomorphisms). This is an increasing function and at each point the derivative is either $\infty$ or positive. In turn, this shows that its inverse,  say $g,$ is differentiable everywhere (thus Lipschitz) and the derivative vanishes at least at all $a_n$. Note that the derivative does not vanish identically as the function is not constant. Choosing conveniently $a_n$ or by translations, we can assume that the derivative at $0$ is positive, thus $g(x)\sim x$ near the origin, but the derivative vanishes wherever the derivative of $f$ is $\infty$ i.e. on a residual set. In conclusion, we cannot claim injectivity of the derivative of a Lipschitz homeomorphism
if it is not bi-Lipschitz (even for  homeomorphism differentiable and $\sim x$). The argument in the above proposition does not work for $h=g, \tau=f$ and $B$ the points where $g'=0$ (in particular it contains all $a_n$). However the question about the differentiability of a Lipschitz mapping on a homeomorphic image of open subset in $\R^n$ is interesting.

%

\end{rem}

\bigskip
\section{Lipschitz type of a given set}

We first analise an example.

\begin{example}\label{components}
Let $\tau$ be the plane $\{ z = 1 \}$ in $\R^3$. 
We take 4 points $P_1$, $P_2$, $P_3$ and $P_4$ on $\tau$ in $\R^3$ 
like in \S 3(2). 
Namely,  
$$
P_1 = (1,1,1), P_2 = (1,-1,1), P_3 = (-1,-1,1) \ {\text and} \ P_4 = (-1,1,1).
$$
Let $\mu$ be the half plane on the plane $\tau$ defined by $y \le - 1$. 
For $t \in I := (- \infty, \infty )$, let $C_t \subset \mu$ be an arc of a circle on $\tau$ 
with centre $(0,t,1)$ of radius $\sqrt{1 + (t + 1)^2}$, connecting $P_2$ and $P_3$. 
Let $Q_t$ be a closed curve on $\tau$ consisting of segments 
$P_1 P_2$, $P_3 P_4$, $P_4 P_1$ and a circular arc $C_t$ for $t \in I$.
We denote by $B_t$, $t \in I$, the real cone with vertex $(0,0,0)$ and base $Q_t$ 
like $B$ in \S 3(2).
Then we can see that for any $t_1 \ne t_2 \in I$, $(\R^3, B_{t_1})$ and $(\R^3, B_{t_2})$ 
are bi-Lipschitz equivalent as germs at $0 \in \R^3$ and $L\mathcal{GD}(B_{t_1}) \ne L\mathcal{GD}(B_{t_2})$. 

Let $J_1 = (-\infty , -\frac{3}{2}] $, $J_2 = (-\frac{3}{2}, -1) $, $J_3 =  \{ -1 \}$ and
$J_4 = (-1 , \infty )$ such that $I = J_1 \cup J_2 \cup J_3 \cup J_4$. 
We put $S_t := L\mathcal{GD}(B_t) \cap \tau$ for $t \in I$. 
Let $m_1 (t)$, $t \in I$, be the number of connected components of the complement 
of $S_t$ in $\tau$. 
Then we can see that $m_1 (t) = 3$, $4$, $2$ and $9$ for $t \in J_i$, 
$i = 1$, $2$, $3$ and $4$, respectively.  
We note that for $t \in J_i$, $i = 1$, $2$, $3$, the complement of  $L\mathcal{GD}(B_t ) \cap \{ z = 0\}$
in $\{ z = 0\}$ is empty, but for $t \in J_4$ it is not empty. 
We next let $m(t)$, $t \in I$, be the number of connected components of the complement 
of $L\mathcal{GD}(B_t )$ in $\R^3$.
Then we can see that $m(t) = 6$, $8$, $4$ and $14$ for $t \in J_i$, $i = 1$, $2$, $3$ and $4$, respectively. 
Therefore we can see that  $(\R^3, L\mathcal{GD}(B_{t_i}))$ and $(\R^3, L\mathcal{GD}(B_{t_j}))$ are not 
topologically equivalent as germs at $0 \in \R^3$ for $t_i \in J_i$ and $t_j \in J_j$ if $i \ne j$.
It follows that they are not bi-Lipschitz equivalent. 
On the other hand, for each $i$, $1 \le i \le 4$,  $(\R^3, L\mathcal{GD}(B_{t_1}))$ and 
$(\R^3, L\mathcal{GD}(B_{t_2}))$ are bi-Lipschitz equivalent as germs at $0 \in \R^3$ if  $t_1$, $t_2 \in J_i$. 
\end{example}

We next introduce a couple of collections of sets, relevant in understanding the Lipschitz type of a given set.
Let
$$
D_A:=\{L\mathcal{GD} (\tau(A)) | \, \tau : (\R^n ,0) \to (\R^n ,0) , \text{ arbitrary bi-Lipschitz mappings}\}. 
$$

\begin{rem} Let $A$ be the cone defined as $A_1^+$ in \S 3(2). 

(1) By the observation in \S 3(2), $D_A$ have elements with different dimensions.

(2) $(\R^3, A)$ and $(\R^3, B_t)$, $t \in I$, are bi-Lipschitz equivalent as germs at $0 \in \R^3$, 
where $B_t$, $t \in I$, are the cones given in Example 6.1. 
Therefore $D_A$ have elements with different numbers of connected componets of the complement
of them in $\R^3$.
\end{rem}

By definition, $D_A=D_{\tau(A)}$ clearly for arbitrary bi-Lipschitz mapping $\tau : (\R^n ,0) \to (\R^n ,0)$.
Moreover a given bi-Lipschitz homeomorphism $\tau : (\R^n ,0) \to (\R^n ,0)$ induces a bijection
$\tau^*:D_A \to D_A$ defined by $\tau^*( L\mathcal{GD}(B))= L\mathcal{GD}(\tau(B))$.

On $D_A$ one can introduce the following equivalence relation:

Two sets $B, C\in D_A$ are {\em equivalent}, denoted by $B\sim C$, if and only if $(\R^n, L\mathcal{GD}(B))$ 
and $(\R^n, L\mathcal{GD}(C))$ are bi-Lipschitz equivalent as germs at $0 \in \R^n$.

\begin{example}
Let $A$ be the cone defined as $A_1^+$ in \S 3(2), and let $B_t$, $t \in I$, be cones given in Example 6.1. 
Then, by the observation in Examplle 6.1, we have $L\mathcal{GD}(B_t) \in D_A$, $t \in I$.
Let $t_1 \in J_i$ and $t_2 \in J_j$, $1 \le i, j \le 4$.
Then $L\mathcal{GD}(B_{t_1}) \sim L\mathcal{GD}(B_{t_2})$ if and only if $i = j$.
\end{example}

Let 
$$
L_A:=\{\tau(A) | \tau : (\R^n ,0) \to (\R^n ,0) , \text{ arbitrary bi-Lipschitz mappings}\} .
$$
On $L_A$ one can introduce the following equivalence relation:

Two sets $B, C\in L_A$ are {\em equivalent}, denoted by $B\sim C$, if and only if $(\R^n, L\mathcal{GD}(B))$ and 
$(\R^n, L\mathcal{GD}(C))$ are bi-Lipschitz equivalent as germs at $0 \in \R^n$.

\begin{example}
Let $A$ and $B_t$, $t \in I$, be the same as in Example 6.3. 
Then we have $B_t \in L_A$, $t \in I$.
Let $t_1 \in J_i$ and $t_2 \in J_j$, $1 \le i, j \le 4$.
Then $B_{t_1} \sim B_{t_2}$ if and only if $i = j$.
\end{example}

For arbitrary sets $B, C$ we have $L_B=L_C \, \text{if and only if} \, C\in L_B \, \text{if and only if} \, B \in L_C$ 
if and only if $(\R^n, B)$ and $(\R^n, C)$ are bi-Lipschitz equivalent as germs at $0 \in \R^n$.

We note that $D_A/\sim = L_A/\sim .$

\begin{rem}
Note that   we can also introduce 
$$
C_A:= \{LD\tau(A) | \, \tau : (\R^n ,0) \to (\R^n ,0) , \text{ arbitrary bi-Lipschitz mappings}\}/\sim .
$$

In this context by $\sim$ we mean the bi-Lipschitz equivalence of the directional cones, 
$LD\tau(A)\sim LD\mu(A)$ if and only if they are bi-Lipschitz equivalent.
In this case we note that if $A$ and $\tau (A)$ satisfy condition $(SSP)$ then their directional cones are bi-Lipschitz equivalent.
\end{rem}

Given $\tau : (\R^n ,0) \to (\R^n ,0) $ a bi-Lipschitz homeomorphism, it is easy to see that $\tau$ induces a bijection 
$\tau^*: L_A \to L_A$ defined by $\tau^*(C)=\tau(C)$. 
There is a natural mapping $\pi : L_A \to D_A$ defined by $\pi (B):= L\mathcal{GD} (B) $ 
and $\pi $ commutes with $\tau^*$ for any bi-Lipschitz homeomorphism $\tau$.

This bijection does not always factorise to $L_A/\sim$. 
However whenever $\tau $ is special,
meaning that $d\tau (L\mathcal{GD} (B) )=L\mathcal{GD} (\tau(B) )$ where $d\tau$ is an Arzela-Ascoli derivative, 
it factorises to the identity mapping.

Similarly such $d\tau$ induces a bijection of $D_A$ defined by $d\tau (L\mathcal{GD} (B) )=L\mathcal{GD} (\tau(B) )$ 
which induces the identity at $D_A/\sim =L_A/\sim$.

\bigskip


\smallskip

\end{document}